\documentclass [11pt]{article}
\newfont{\bbb} {msbm10}
\newcommand{\Bbb}[1]{\mbox{\bbb#1}}
\newcommand{\R}{\Bbb{R}}
\newcommand{\bS}{\Bbb{S}}

\newcommand{\Z}{\Bbb{Z}}
\newcommand{\Q}{\Bbb{Q}}
\newcommand{\D}{\Bbb{D}}
\newcommand{\cT}{{\cal{T}}}
\newcommand{\sm}{\setminus}
\newcommand{\sbs}{\subset}
\newcommand{\ra}{\rightarrow}

\newcommand{\cP}{{\cal{P}}}

\newcommand{\p}{\partial}

\newcommand{\cM}{{\cal{M}}}
\newcommand{\met}{{\cal{MET}}}
\newcommand{\cD}{{\cal{D}}}

\newcommand{\cG}{{\cal{G}}}
\newcommand{\cA}{{\cal{A}}}

\newcommand{\td}{{\widetilde{DIFF}}}
\usepackage[all]{xy}

\begin{document}

\title{The Teichm\"{u}ller Space of Pinched Negatively Curved Metrics on a Hyperbolic
Manifold is not Contractible}
\author{F. T. Farrell and P. Ontaneda\thanks{The first author was
partially supported by a NSF grant.
The second author was supported in part 
by research grants from CAPES(Brazil) and NSF.}}
\date{}

\maketitle

\begin{abstract} For a smooth manifold $M$ we define the Teichm\"uller space $\cT(M)$
of all Riemannian metrics on $M$ and the Teichm\"uller space $\cT^\epsilon(M)$
of $\epsilon$-pinched negatively curved metrics on $M$, where $0\leq\epsilon\leq\infty$.
We prove that if $M$ is hyperbolic the natural inclusion
$\cT^\epsilon(M)\hookrightarrow\cT(M)$ is, in general, not homotopically trivial. In particular,
$\cT^\epsilon(M)$ is, in general, not contractible. 
\end{abstract}
\vspace{.3in}

\noindent {\bf \large  Section 0. Introduction.}\\

Let $M$ be a closed smooth manifold. We denote by $\met (M)$ the space of all smooth Riemannian metrics
on $M$ with the smooth topology. Note that the space $\met ( M ) $ is contractible.
We also denote by  $DIFF(M)$ the group of all smooth self-diffeomorphisms of $M$.
We have that $DIFF(M)$ acts on $\met (M)$ pulling-back metrics:
$\phi g =(\phi^{-1})^*g=\phi_*g$, for $g\in\met(M)$ and $\phi\in DIFF(M)$,
that is, $\phi g $ is the metric such that $\phi :(M,g)\ra (M,\phi g)$ is an isometry.\\

Let $\R^+$ be the set of positive real numbers which we consider as a group with
multiplication. Denote by $\cD(M)$ the group $\R^+\times DIFF(M)$.
The group $\cD(M)$ acts on $\met(M)$ by scaling and pulling-back metrics:
$(\lambda ,\phi ) g =\lambda (\phi^{-1})^*g=\lambda\phi_*g$, for $g\in\met(M)$ and $(\lambda, \phi)\in\cD(M)$.
The quotient space $\cM(M)=\met(M)/\cD(M)$ is called the moduli space of metrics on $M$.
It is sometimes said that a geometric property is a property that is invariant by isometries,
that is, by the action of $DIFF(M)$. Hence if two Riemannian metrics represent the
same element in $\cM(M)$ then they posses the same geometric properties.
Clearly, the study of the moduli space of metrics is of fundamental importance not just in
geometry but in other areas of mathematics as well. (See, for instance, \cite{B} Ch. 4)\\

It is also interesting to consider subspaces of $\cM(M)$ that represent some geometric property.
One obvious choice is to consider metrics with constant curvature. For instance,
let $M_g$ be an orientable two-dimensional manifold of genus $g>1$. 
Consider the moduli space of all hyperbolic metrics on $M_g$, that is, the subspace of $\cM(M_g)$
formed by elements that are represented by Riemannian metrics of constant sectional curvature equal to -1.
The moduli space of all hyperbolic metrics
is the quotient of another well known space: the Teichm\"uller space of $M_g$. This space is a subspace of
the quotient of $\met(M_g)$ by the subgroup of $DIFF(M_g)$ formed by all smooth self-diffeomorphisms of $M_g$
which are homotopic to the identity; namely, it is the subspace represented by hyperbolic metrics.
Then the moduli space is the quotient of the Teichm\"uller space by the
action of $Out\, (\pi_1(M_g))$, the group of outer automorphisms of the fundamental group of $M_g$.\\

We want to generalize the definition of the Teichm\"uller space to higher dimensions.
The obvious choice for a definition would be the quotient of the space of all hyperbolic metrics
by the action of the group of all smooth self-diffeomorphisms which are homotopic to the identity.
But Mostow's Rigidity Theorem implies that, in dimensions $\geq 3$, this space contains (at most) one point.\\

Let us go back to dimension two for a moment. Recall that uniformization techniques (see \cite{EE}, or,
more recently, Hamilton's Ricci flow \cite{H}) show that every Riemannian metric on $M_g$, $g>1$,
can be canonically deformed to a hyperbolic metric. Moreover, Hamilton's Ricci flow \cite{H}
shows that every negatively curved metric on $M_g$, $g>1$,
can be canonically deformed (through negatively curved metrics) to a hyperbolic metric. Hence the space of all 
hyperbolic metrics on $M_g$ is
canonically a deformation retract of the space of all negatively curved Riemannian metrics on $M_g$.
This deformation commutes with the action of $DIFF(M_g)$ (this is true at least for the Ricci flow), therefore
the Teichm\"uller space of $M_g$ is canonically a deformation retract of the space which is the quotient
of all negatively curved Riemannian metrics on $M_g$ by the action of 
the group of all smooth self-diffeomorphisms which are homotopic to the identity. Also, instead of considering
the space of all negatively curved metrics we can consider the space of all {\it pinched }
negatively curved metrics, or for that matter, the space of all Riemannian metrics.
These are the concepts that we will generalize. Next, we give detailed definitions
and introduce some notation.\\
 
As before, let $M$ be a closed smooth manifold. We denote by $DIFF\, _0(M)$ the subgroup of $DIFF(M)$ of
all smooth diffeomorphisms of $M$ which are homotopic to the identity $1_M$. 
Also, denote by $\cD_0(M)$ the group $\R^+\times DIFF_0(M)$. 
We call the quotient space $\cT(M)=\met(M)/\cD_0(M)$ the Teichm\"{u}ller space of metrics on $M$.\\

Given $0\leq \epsilon\leq \infty$ let $\met^\epsilon (M)$ denote the space of all 
$\epsilon$-pinched negatively curved Riemannian metrics on $M$, that is, $\met^\epsilon (M)$ is the
space of all negatively curved Riemannian metrics $g$ on $M$ for which
$$\frac{sup{\mbox{-}}sec\, g}
{ inf{\mbox{-}}sec\, g}\,\, \leq\,\, 1+\epsilon$$

\noindent where $sup{\mbox{-}}sec \, g=sup\{ |g{\mbox{-sectional curvature of }}P_x|,\,\,\, x\in M, P_x\, 
{\mbox{two-plane in}}\, T_x M\}$
and $inf{\mbox{-}}sec \, g=inf\{ |g{\mbox{-sectional curvature of }}P_x|,\,\,\, x\in M, P_x\, 
{\mbox{two-plane in}}\, T_x M\}$.
Therefore $g\in \met^\epsilon (M)$ if and only if there is a positive real number $\lambda$ such that
$\lambda g$ has all its sectional curvatures in the interval $[-(1+\epsilon ), -1]$.
Note that a 0-pinched metric is a metric of constant negative sectional curvature and an
$\infty$-pinched metric is just a negatively curved Riemannian metric.\\

The quotient space $\cM^\epsilon (M)=\met^\epsilon (M)/\cD(M)$ 
is called the moduli space of $\epsilon$-pinched negatively curved metrics on $M$. Also, 
$\cT^\epsilon (M)=\met^\epsilon (M)/\cD_0(M)$ is called the Teichm\"{u}ller space of $\epsilon$-pinched 
negatively curved metrics on $M$. In particular, $\cT^\infty (M)$ is the Teichm\"uller space of all
negatively curved metrics on $M$.
Note that the inclusions $\met^\epsilon (M)\hookrightarrow \met (M)$ induce inclusions
$\cT^\epsilon (M)\hookrightarrow \cT (M)$.
Also note that, for $\delta\geq\epsilon$,  these inclusions factor as follows:
$\met^\epsilon (M)\hookrightarrow \met^\delta (M)\hookrightarrow\met (M)$ and
$\cT^\epsilon (M)\hookrightarrow \cT^\delta (M)\hookrightarrow\cT (M)$.\\

\noindent {\bf Remarks.}\\

\noindent {\bf A.} 
If $M_g$ is an orientable two-dimensional manifold of genus $g>1$, then the original Teichm\"uller space of
$M_g$ is denoted (in our notation) by $\cT^0(M_g)$ and
$\cT^0(M_g)$ is homeomorphic to $\R^{6g-6}$ (see \cite{EL2}).
Hence $\cT^0(M_g)$ is contractible. By the uniformization techniques mentioned above (\cite{EE}, \cite{H}), it 
follows that $\cT^\epsilon(M_g)$, $\cT^\infty (M_g)$, $\cT(M_g)$ are all contractible. 
(This is also true for non-orientable surfaces of Euler characterisitc $<0$.)\\

\noindent {\bf B.}
Let $M$ be a closed hyperbolic manifold. If 
$dim\, M\geq 3$, Mostow's Rigidity Theorem implies that $\cT^0(M)=*\,$; ie. $\cT^0(M)$ contains exactly
one point. Therefore $\met^0(M)=\cD_0(M)$. It also follows (see A above) that
$\cT^0(M)$ is contractible when $dim\, M\geq 2$.\\

In dimensions two and three it is known that $\cD_0(M)$ (and hence $\met^0(M)$) is
contractible. (This is due to Earle and Eells \cite{EE} in dimension two and to Gabai \cite{G} in dimension three.)
This is certainly false in dimensions $\geq 11$,
because $\pi_0(\cD_0(M))$ is not finitely generated (see \cite{FJ3}, Cor.10.16 and 10.28), and it is reasonable
to conjecture that $\cD_0(M)$ is also not contractible for dimension $n$, $5\leq n\leq 10$.\\

\noindent {\bf C.}
Let $M$ be a hyperbolic manifold. Then the action of $\cD_0(M)$ on $\met(M)$ is free
(see Lemma 1.1). Since $\met(M)$ is contractible 
by Ebin's Slice Theorem \cite{Eb}
we have that $\cD_0(M)\ra\met (M)\ra\cT (M)$
is a principal $\cD_0(M)$-bundle and $\cT(M)$ is the classifying space $B\cD_0(M)$ of
$\cD_0(M)$. \\

Therefore, if $M$ is a closed hyperbolic manifold then $\met^\epsilon (M)$ interpolates between $\met^0(M)$
(which is equal to $\cD_0(M)$) and $\met (M)$ (which is contractible). Likewise 
$\cT^\epsilon (M)$ interpolates between $\cT(M)$ (which is equal to $ B\cD_0(M)$)
and $\cT^0 (M)$ (which is contractible). Schematically, we have the following diagram:\\

$$\begin{array}{ccccccc}
\met^0(M)&\hookrightarrow & \met^\epsilon (M)&\hookrightarrow &\met^\infty (M)&\hookrightarrow &\met(M)\\
\downarrow & &\downarrow &&\downarrow &&\downarrow \\
\cT^0(M)&\hookrightarrow & \cT^\epsilon (M)&\hookrightarrow &\cT^\infty (M)&\hookrightarrow &\cT(M)
\end{array}$$\\

\noindent All vertical arrows represent quotient maps by the action of the group $\cD_0(M)$.\\

The main result of this paper states that for a hyperbolic manifold the 
last two horizontal arrows of the lower row of the diagram above are not in general homotopic to a constant
map. In particular $\cT^\epsilon$, $0\leq\epsilon\leq\infty$, is in general not contractible.
More specifically, we prove that under certain conditions on the dimension $n$ of the hyperbolic manifold
$M$, the manifold $M$ has a finite cover $N$ (which depends on $\epsilon$) such that 
$\pi_k (\cT^\epsilon (N))\ra\pi_k(\cT (N))$ is non-zero. In particular, $\cT^\epsilon (N)$ is not contractible.
The requirements on the dimension $n$ are implied by one of the following conditions:
$n$ is larger than some constant $n_0(4)$ or $n$ is larger than 5 but in this last case we need
that $\Theta_{n+1}\neq 0$, where $\Theta_{\ell }$ denotes the group of homotopy spheres of dimension $\ell$.
Here is a more detailed statement of our main result:\\

\noindent {\bf Theorem 1.} {\it For every integer $k_0\geq 1$ there is an integer $n_0=n_0(k_0)$ such that
the following holds. Given $\epsilon >0$ and a closed real hyperbolic $n$-manifold $M$ with  $n \geq n_0$,
there is a finite sheeted cover $N$ of $M$ such that, for every $1\leq k\leq k_0$ with $n+k\equiv 3$ $mod\,\, 4$, the map
$\pi_k (\cT^\epsilon (N))\ra\pi_k(\cT (N))$, induced by the inclusion $\cT^\epsilon (N)\hookrightarrow \cT (N)$,
is non-zero.
Consequently $\pi _k (\cT^\epsilon (N))\neq 0$. In particular, 
$\cT^\delta (N)$ is not contractible, for every $\delta$ such that $\epsilon\leq\delta\leq\infty$
(provided $k_0\geq 4$).}\\

Here (and in the Corollary below) we consider the given hyperbolic metric as the basepoint for
$\cT (N)$, $\cT^\epsilon (N)$. \\

For $k_0=1$ we will show that we can take $n_0(1)=6$, and that we can drop the condition $n+k\equiv 2$ $mod\,\, 4$.
Hence we obtain the following corollary to (the proof of) Theorem 1.\\

\noindent {\bf Corollary.} {\it Let $M$ be a closed real hyperbolic manifold of dimension $n$, $n\geq 6$.
Assume that $\Theta_{n+1}\neq 0$.
Then for every $\epsilon >0$ there is a finite sheeted cover $N$ of $M$ such that 
$\pi _1 (\cT^\epsilon (N))\neq 0$. Therefore
$\cT^\epsilon (N)$ is not contractible.}\\

Recall that an $n$-dimensional $\pi$ manifold is a manifold that embeds in $\R^{2n+2}$ with
trivial normal bundle. Every real hyperbolic manifold has a finite sheeted cover that is
a $\pi$ manifold (see \cite{Su}, p.553). We have the following addition to the statements of Theorem 1 and
the Corollary.\\

\noindent {\bf Addendum.} {\it We can choose $N=M$ in the statements of
Theorem 1 and the Corollary, provided $M$ is a $\pi$-manifold and the radius of injectivity of $M$ at some point is sufficiently
large (how large depending only on the dimension of $M$).}\\

We now make some comments on Theorem 1 and  the diagram above.\\

\noindent {\bf 1.} Since $\met(M)$ is contractible, Theorem 1 implies that, for a general
hyperbolic manifold $M$,
the map $\pi_k(\met^\epsilon(M))\ra\pi_k(\cT^\epsilon(M))$, induced by the second vertical
arrow of the diagram, is not onto for some $k$.\\

\noindent {\bf 2.} By remark A, 
the lower row of the diagram above is homotopically trivial in dimension 2.
In dimension 3 one could ask the same:
is the lower row of the diagram above homotopically trivial in dimension 3?
In view of a result of Gabai (see \cite{G}), this is equivalent to asking: is $\cT(M)^\infty$
contractible?\\

\noindent {\bf 3.} Let $M$ be a hyperbolic manifold. Consider the upper row of the diagram.
It follows from a result of Ye on the Ricci flow (see \cite{Ye}) that, provided the dimension of $M$ is even,
there is an $\epsilon_0=\epsilon_0(M)>0$ such that for all $\epsilon\leq\epsilon_0$ the inclusion map
$\met^\epsilon\ra\met^\infty$ is $\cD_0(M)$-equivariantly homotopic to a retraction 
$\met^\epsilon\ra\met^0(M)\sbs\met^\infty$. This has the following consequences.
First the retraction above descends to a retraction $\cT^\epsilon(M)\ra\cT^0(M)$, hence the inclusion map
$\cT^\epsilon(M)\ra\cT^\infty(M)$ is homotopic to a constant map (provided $\epsilon\leq\epsilon (M)$),
and hence induces the zero homomorphism $\pi_k(\cT^\epsilon(M))\ra\pi_k(\cT(M))$ for all $k$.
Second, the inclusion map $\cD_0(M)=\met^0(M)\ra\met^\epsilon(M)$ induces monomorphisms
$\pi_k(\cD_0(M))=\pi_k(\met^0(M))\ra\pi_k(\met^\epsilon(M))$, provided $\epsilon\leq\epsilon(M)$.
Theorem 1 then shows that in many cases
$\epsilon_0(M)<\infty$. \\

\noindent {\bf 4.} We recall now an open problem posed by K. Burns and A. Katok (\cite{BK}, Question 7.1)
about hyperbolic manifolds $M$. Is $\met^\infty(M)$ path connected?
More generally one could ask if $\met^\infty(M)$ is contractible. Equivalently, is $\cT^\infty (M)\ra\cT(M)=B\cD_0(M)$
a homotopy equivalence? Before Theorem 1, it was conceivable that the opposite extreme in the interpolation between
$\cT^0(M)$ and $\cT(M)$ could be true, i.e. $\cT^\infty(M)$ is always contractible, or equivalently
$*=\cT^0(M)\ra\cT^\infty(M)$ is always a homotopy equivalence.\\

\noindent {\bf 5.} Let $M$ be a hyperbolic manifold. Since $DIFF(M)/DIFF_0(M)\cong Out (\pi_1(M))$
we have that $\cM(M)\cong\cT(M)/Out(\pi_1(M))$ or, in general, $\cM^\epsilon (M)\cong\cT^\epsilon (M)/Out(\pi_1(M))$.
Note that $Out (\pi_1(M))$ is a finite group, provided $dim\, M\geq 3$. We do not know whether our results descend to the moduli
spaces.\\

\noindent {\bf 6.} Let $M$ be a hyperbolic manifold.
We can consider the quotients of $\met(M)$ and $\met^\epsilon (M)$ by
$DIFF^0(M)$, the connected component of the identity
$1_M$ in $DIFF(M)$, instead of by the larger group $DIFF_0(M)$. Since the quotient
group $DIFF_0(M)/DIFF^0(M)$ is discrete, it can be easily checked from the proof of our results that
the statement of Theorem 1 also holds for the inclusion of the quotients:
$\met^\epsilon (M)/DIFF^0(M)\ra \met (M)/DIFF^0(M)$, with the strengthened restriction
``$\, 2\leq k\leq k_0$'' and proviso ``(provided $k_0\geq 5$)''.\\


Theorem 1 follows from the more technical Theorems 2 and 3 below. To state these results we need some notation.
Write $\cG\sbs DIFF_0\, (\bS^{n-1}\times [1,2], \p )$ for the group of all smooth isotopies $\varphi$ of the $(n-1)$-dimensional
sphere $\bS^{n-1}$ which are the identity near 1,2, and are homotopic to the identity
by a homotopy that is constant near 1,2.
That is, $\varphi : \bS^{n-1}\times [1,2]\ra\bS^{n-1}\times [1,2]$,
$\varphi (x,s)=(y,s)$ and $\varphi (x,s)=(x,s)$ for $s$ near 1,2 and $\varphi$ is homotopic to the identity
by a homotopy $H$ such that $H_t(x,s)=(x,s)$ for $s$ near 1,2, and all $t\in[0,1]$. 
Note that $\cG$ depends only on the dimension $(n-1)$ of the sphere. If we need to express this
dependency explicitly we shall write $\cG_{n}$.\\

Let $N$ be a real hyperbolic manifold of dimension $n$ and let $B$ be a closed geodesic ball of radius $2r$, 
centered at some point $p\in N$ (hence, the radius of injectivity of $N$ at $p$ is larger that $2r$). 
Identify $B\sm\{ p\}$ with  $\bS^{n-1}\times (0,2r]$, where the lines $t\mapsto (x,t)$ are the
speed-one geodesics emanating from $p$. Now, every element in $DIFF_0(\bS^{n-1}\times [1,2],\p )$ gives rise to an element in
$DIFF_0(N)$ by identifying $\bS^{n-1}\times [1,2]$ with $\bS^{n-1}\times [r,2r]$. That is, we have a map
$\Lambda =\Lambda (N, p, r):DIFF_0(\bS^{n-1}\times [1,2],\p )\ra DIFF_0(N)$, defined by\\

$$\Lambda \varphi (p)=\left\{\begin{array}{cll} p &  &p\notin (\bS^{n-1}\times [r,2r])\sbs N \\ \\
(\varphi_{t/r}(x), t) && p=(x,t)\in (\bS^{n-1}\times [r,2r])\sbs N\end{array}\right.$$\\

\noindent where $\varphi\in DIFF_0(\bS^{n-1}\times [1,2],\p )$ and $\varphi (x,s)=(\varphi_s(x),s)$.
We will denote the restriction $\Lambda|_{\cG}$ by the same symbol $\Lambda$.\\

\noindent {\bf Remarks.}  

\noindent {\bf 1.}    
A technical point here. Rigorously, for $\Lambda$ to be well defined (i.e. for $\Lambda\varphi$ to be smooth)
we will assume that every element in $DIFF_0(\bS^{n-1}\times [1,2],\p )$ is the identity near 
$\bS^{n-1}\times \{ 1\}$ and $\bS^{n-1}\times \{ 2\}$.
This does not cause problems since standard extension methods (along collars) show that the 
inclusion of the space of all elements in $DIFF_0(\bS^{n-1}\times [1,2],\p )$, with the properties above, into 
$DIFF_0(\bS^{n-1}\times [1,2],\p )$ is a homotopy equivalence.\\

\noindent {\bf 2.} Another technical point. The map $\Lambda$  depends also on the identification
between $B\sm\{ p\}$ and  $\bS^{n-1}\times (0,2r]$ used above. This identification is uniquely
determined if an orthonormal basis $\cal B$ of $T_pM$ is specified. Hence we should write
$\Lambda =\Lambda (N, p, r, {\cal{B}})$. To alleviate the notation we still
write $\Lambda =\Lambda (N, p, r)$ since the choice of $\cal B$ is not essential.
Note also that such a map $\Lambda$ can be defined whenever the radius of injectivity of $N$ at $p$
is larger than $2r$.\\

\noindent {\bf Theorem 2.} {\it Given $\epsilon >0$ and a compact subset $K\sbs\cG$ there is a 
real number $r>0$
such that the following holds. Let $(N,g^0)$ be a closed real hyperbolic manifold and $p\in N$ with
radius of injectivity at $p$ larger that $3r$. Then the map $K\ra \met^\epsilon(N)$ given by
$\phi\mapsto (\Lambda \phi )g^0$
is contractible; i.e. homotopic to a constant map. Here $\Lambda=\Lambda(N,p,r)$.}\\

\noindent {\bf Theorem 3.} {\it For every integer $k\geq 0$ there is an integer $n_1=n_1(k)$ 
and elements $\alpha_{k,n}\in\pi_k(\cG_n)$, $n\geq n_1$, such that
the following holds. If $N$ is a closed real hyperbolic $n$-manifold, with $n\geq n_1$, $n+k\equiv 2$ $mod\,\, 4$,
which is a $\pi$-manifold then  $\Lambda_\# (\alpha_{k,n})\neq 0
\in\pi_k(DIFF_0(N))$. Here $\Lambda=\Lambda (N,p,r)$, where $p\in N$ is any point and $2r$ is 
less than the injectivity radius of $N$ at $p$.}\\

For the proof of the Corollary we will need the following Addendum to Theorem 3.\\

\noindent {\bf Addendum.} {\it  For $k=0$ we can choose $n_1=n_1(1)=6$ and we can drop the condition
$n+k\equiv 2$ $mod\,\, 4$, provided $\Theta_{n+1}\neq 0$}.\\

Note that, while Theorem 2 is a geometric result, Theorem 3 is purely topological: the point $p$ and the number
$r$ are arbitrary, and the only restriction is that $2r$ is less than the injectivity radius of $N$ at $p$.\\

In Section 1 we deduce Theorem 1 and the Corollary from Theorems 2 and 3. In Section 2 we prove Theorem 2 and 
in Section 3 we prove Theorem 3 together with its Addendum.\\

We are grateful for the referee's useful comments and suggestions.
\vspace{.5in}

\noindent {\bf \large  Section 1. Proof of Theorem 1 and the Corollary.}\\

Here we prove Theorem 1 and the Corollary assuming Theorems 2 and 3. First we give a Lemma and some remarks.\\

Recall that $DIFF_0(P)$ and $\cD_0(P)$ act on $\met(P)$, for any closed smooth manifold $P$.\\

\noindent {\bf Lemma 1.1.} {\it If $P$ is aspherical and the center of $\pi_1 P$ is trivial then
the action of $DIFF_0(P)$ and $\cD_0(P)$ on $\met (P)$ is free.}\\

\noindent {\bf Proof.} Let $g\in\met(P)$. Note that the isotropy group $H=\{ \phi\in DIFF_0(P)\,,$
$\, \phi g =g\}$ of the action of $DIFF_0(P)$ at $g$ is $ISO_0(M,g)$, the group of all isometries of the Riemannian
manifold $(M,g)$ that are homotopic to the identity. Hence this isotropy group $H$ is compact. Let
$\gamma :DIFF(P)\ra Out\,(\pi_1\, P)$ be the homomorphism induced by $\phi\mapsto \phi_*$. Borel-Conner-Raymond
showed (see \cite{CR}, p. 43) that under the assumptions above, $\gamma$ restricted to compact subgroups is monic.
But $\gamma (H)$ is trivial, since every element in $DIFF_0(P)$ is, by definition, homotopic to the
identity. It follows that $H$ is trivial. Hence the action of $DIFF_0(P)$ is free.
Therefore the action of $\cD_0(P)$ is also free. This proves the Lemma.\\

It follows from the Lemma that there is a fibration $\cD_0(P)\ra\met(P)\ra\cT(P)$, and, since $\met(P)$ is contractible,
we have that $\pi_{k-1}(DIFF_0(P))\cong\pi_{k-1}(\cD_0(P))\cong\pi_k(\cT(P))$.
We can give an explicit isomorphism between $\pi_{k-1}(DIFF_0(P))$ and $\pi_k(\cT(P))$. Let $g^0$ be any metric on $P$.
Let $\beta :\bS^{k-1}\ra DIFF_0(P)$ be an element in $\pi_{k-1}(DIFF_0(P))$. 
Define $\beta':\bS^{k-1}\ra\met(P)$ by $\beta' (u)=\beta(u)g^0$, $u\in\bS^{k-1}$.
Since $\met(P)$ is contractible
we can extend $\beta'$ to a map $\bar\beta$ defined on the whole disc $\D^k$. Then the isomorphism
is given by $\beta\mapsto p\bar\beta$, where $p:\met(P)\ra\cT(P)$ is the quotient map that assigns to each metric
its $\cD_0(P)$ orbit.\\

\noindent {\bf Proof of Theorem 1 assuming Theorems 2 and 3.} 
Let $\epsilon >0$, $k_0>0$. 
Let  $n_0=n_0(k_0)=max\{ n_1(k-1)\, ,\,\, 1\leq k\leq k_0\}$, where $n_1(k-1)$ is as in Theorem 3. Let $n$
be such that $n\geq n_0$.
Define the compact subset $K$ of $\cG$ by
$$K=\{ \alpha_{k-1}(u)\, ,\  u\in\bS^{k-1}\,\, ,\,\, 1\leq k\leq k_0\, , \,\,  n+(k-1)\equiv 2 \,\,mod\,\, 4
\}$$

\noindent where the $\alpha_{k-1}=\alpha_{(k-1),n}$ are explicit representatives of the elements given in Theorem 3.
Note that $K=\cup_{k} \,\,image(\alpha_{k-1})$ where $1\leq k\leq k_0$ and $ n+(k-1)\equiv 2 \,\,mod\,\, 4$.
Let $r$ be as in Theorem 2, for $\epsilon$ and $K$ as above. \\

Now, let $M$ be a closed hyperbolic manifold of dimension $n$. By taking successive finite sheeted
covers we can find a finite sheeted cover $N$ of $M$ such that:
\begin{enumerate}
\item[{1.}] $N$ is a $\pi$ manifold (see \cite{Su}, p. 553).

\item[{2.}] $N$ has a point $p$ with injectivity radius larger than $3r$. (Recall that $\pi_1(M)$ is
residually finite, see \cite{Mag}.)
\end{enumerate}

Let $g^0$ be the hyperbolic metric of $N$ pulled back from $M$.
Write $\Lambda=\Lambda(N,p,r)$. Define $\beta_{k-1}=\Lambda\alpha_{k-1}\in\pi_{k-1}(DIFF_0(N))$. By Theorem 3 all
$\beta_{k-1}$ are non-zero. Define also $\beta_{k-1}':\bS^{k-1}\ra \met^\epsilon (N)$ by $\beta_{k-1}' (u)=\beta_{k-1}(u)g^0$, 
$u\in\bS^{k-1}$.
By Theorem 2 we can extend each $\beta'=\beta'_{k-1}$ to the whole disc $\D^k$, obtaining maps
$\bar\beta :\D^k\ra\met^\epsilon(N)$.\\

Recall that $p:\met(P)\ra\cT(P)$ is the quotient map that assigns to each metric
its $\cD_0(P)$ orbit (see the comments following the proof of Lemma 1.1).
Since $p\bar\beta ( \bS^{k-1}) $ contains exactly one point (this point is $p(g^0)$) we have that
$p\bar\beta$ determines an element in $\pi_k(\cT(N))$. Also, since
$\beta_{k-1}\neq 0$ we have that $p\bar\beta\neq 0\in \pi_k(\cT(N))$.
But $image\, (\bar\beta)\sbs\met^\epsilon(N)$, hence we have that $p\bar\beta$ is in the image
of the map $\pi_k (\cT^\epsilon (N))\ra\pi_k(\cT (N))$ induced by the inclusion $\cT^\epsilon (M)\hookrightarrow \cT (M)$.
This proves Theorem 1 assuming Theorems 2 and 3.\\

The proof of the Corollary is similar, just use the Addendum to Theorem 3. 
The proof of the Addendum (to Theorem 1 and the Corollary) is also similar.
\vspace{.4in}

\noindent {\bf \large  Section 2. Proof of Theorem 2.}\\

First we introduce some notation and give a Lemma.
Denote by $\cG'\sbs DIFF_0\, (\bS^{n-1}\times [1,2])$ the group of all smooth isotopies $\varphi$ of the $(n-1)$-dimensional
sphere $\bS^{n-1}$ which are the identity near 1, and constant near 2.
That is, $\varphi : \bS^{n-1}\times [1,2]\ra\bS^{n-1}\times [1,2]$,
$\varphi (x,s)=(y,s)$, $y$ does not depend on $s$, for $s$ near 2 and $\varphi (x,s)=(x,s)$ for $s$ near 1. 
Note that $\cG'$ depends only on the dimension $(n-1)$ of the sphere.
We have an inclusion $\cG\hookrightarrow\cG'$.\\

\noindent {\bf Lemma 2.2.} {\it  $\cG'$ is contractible.}\\

\noindent {\bf Proof.} Recall that the space of all isotopies is homeomorphic to the space of smooth paths of
vector fields $V_s$, $s\in[1,2]$ on the sphere (or any closed manifold). This correspondence is given explicitly in the 
following way. An isotopy $\varphi$ corresponds to the smooth path of vector fields $V_s$, where
$V_s(x)=\frac{d}{ds}\varphi_s(x)|_{s}$. Here $\varphi (x,s)=(\varphi_s(x),s)$. Conversely, given
a smooth path of vector fields $V_s$, $s\in[1,2]$, we can integrate it and obtain the flow 
$\varphi_s$ of $V_s$. Then $V_s$ corresponds to the isotopy $\varphi (x,s)=(\varphi_s(x),s)$.
But every vector field (or path of vector fields) can be deformed to the zero vector field
by homotheties: $(\mu, V)\mapsto (1-\mu)V$, $\mu\in[0,1]$, is a homotopy of $V$ to the zero vector field.
Integrating this homotopy (for each $\mu$) we obtain a homotopy from the identity $id_{\cG'}$ to the
constant map $\cG'\ra\{ id_{\bS^{n-1}\times [1,2]}\}\sbs\cG'$. This proves the Lemma.\\

Note that the homotopy given in the proof of the Lemma fixes the identity $id_{\bS^{n-1}\times [1,2]}$.
Note also that the homotopy does not necessarily leave $\cG$ invariant.\\

{\bf Proof of Theorem 2.} We first prove the Theorem for the case in which $K$ has exactly one element.
Fix $\epsilon>0$ and $\varphi\in\cG$. Let $(N,g^0)$ be a real hyperbolic manifold and $p\in N$ with
injectivity radius (at $p$) larger than $3r$.
We will construct a deformation of $(\Lambda\varphi ) g^0=(\Lambda\varphi)_*g^0$ to $g^0$, through metrics in $\met^\epsilon(N)$,
assuming that $r$ is large enough. Here $\Lambda=\Lambda(N,p,r)$.\\

Denote by $B\sbs N$ the closed geodesic ball centered at $p$ of radius $3r$. As before, we identify
$B\sm\{ p\}$ with $\bS^{n-1}\times (0,3r]$. In fact, this identification can be done isometrically:
$B\sm\{ p\}$ with metric $g^0$ is isometric to $\bS^{n-1}\times (0,3r]$ with metric 
$sinh^2(t)h+dt^2$, where $h$ is the Riemannian metric on the sphere $\bS^{n-1}$ with constant
curvature equal to 1. In view of this identification we write then $g^0(x,t)=sinh^2(t)h(x)+dt^2$.\\

Write $\phi=\Lambda\varphi$. Also
write $g^1=\phi  g^0$. The metric $g^1$ on $B\sm\{ p\}$ is given by:

$$g^1(x,t)=\left\{\begin{array}{lll} g^0(x,t) &  &t\notin [r,2r] \\ \\
\phi_*g^0\, (x,t) && t\in  [r,2r]\end{array}\right.$$\\

By the above Lemma 2.2 we have a path of isotopies $\varphi^\mu\in\cG'$, $\mu\in[0,1]$, with 
$\varphi^0=\varphi$ and $\varphi^1=id_{\bS^{n-1}\times [1,2]}$. Write $\theta^\mu=\varphi^\mu_2$ for the
final map of the isotopy $\varphi^\mu$, that is, $\varphi^\mu(x,2)=(\theta^\mu(x),2)$.
Then $\theta^0=\theta^1=id_{\bS^{n-1}}$.
Define $\phi^\mu:\bS^{n-1}\times[r,2r]\ra \bS^{n-1}\times[r,2r]$ by rescaling $\varphi$ to the interval
$[r,2r]$, that is, $\phi^\mu(x,t)=(\varphi^\mu_{t/r}(x),t)$. \\

Let $\delta :[2,3]\ra[0,1]$ be smooth with $\delta(2)=1$, $\delta(3)=0$ and $\delta$ is constant near 2 and 3.
We now define a path of metrics $g^\mu$ on $B\sm\{ p\}=\bS^{n-1}\times (0,3r]$:\\

$$g^\mu(x,t)=\left\{\begin{array}{lll} g^0(x,t) &  &t\in (0,r] \\ \\
(\phi^\mu)_*g^0\, (x,t) && t\in  [r,2r]\\ \\
sinh^2(t)\Bigl[ \delta(\frac{t}{r})(\theta^\mu)_*h(x)+ \{1-\delta(\frac{t}{r})\} h(x))\Bigr] +dt^2 &&t\in  [2r,3r]
\end{array}\right.$$\\

Since $\delta$ and all isotopies we used are constant near the endpoints of their intervals of
definitions, it is straightforward to show that $g^\mu$ is a smooth metric on $B\sm\{p\}$ and that
$g^\mu$ joins $g^1$ to $g^0$. Moreover,  $g^\mu(x,t)=g^0(x,t)$, for $t$ near 0 and 3. Hence we can extend
$g^\mu $ to the whole manifold $N$ by defining $g^\mu(q)=g^0(q)$, for $q=p$ or $q\notin B$.\\

\noindent {\bf Claim.} {\it $g^\mu\in\met^\epsilon(N)$, provided $r$ is large enough 
(depending only on $\varphi^\mu$ and $\delta$).}\\

\noindent {\bf Proof.} Since $g^\mu(x,t)$ is equal to $g^0(x,t)$, for $t\in(0,r]$, hence
$g^\mu(x,t)$ is hyperbolic for $t\in(0,r]$. Also, $g^\mu(x,t)$ is the push-forward (by $\phi^\mu$) of the hyperbolic
metric $g^0$, for $t\in[r,2r]$; hence $g^\mu(x,t)$ is hyperbolic, for $t\in[r,2r]$. For 
$t\in[2r,3r]$, the metric $g^\mu(x,t)$ is similar to the ones constructed in
\cite{FJ1}, Section 3,  or \cite{O}, Th. 3.1. It can be checked from those references that the
sectional curvatures of $g^\mu$ are $\epsilon$ close to -1, 
provided $r$ is large enough. How large we need $r$ to be depends only on 
the partial derivatives (up to order two) of $\varphi^\mu$ and $\delta$. This proves the claim
and Theorem 2 for the case in which the compact set $K$ has exactly one element.\\

For the general case just note that since $K$ is compact so is the set $\bar K =\{ \varphi^\mu\, ,\, \varphi\in K\}$,
where $\varphi^\mu$ denotes the canonical deformation of an element $\varphi\in K$ to the identity
(given by Lemma 2.2 above). Then all partial derivatives (up to order two) of all elements in $\bar K$ are
bounded. Therefore there is a real number $r$ for which the argument used in the claim above works
for all $\varphi\in K$. This proves Theorem 2.
\vspace{.4in}

\noindent {\bf \large  Section 3. Proof of Theorem 3.}\\

In this section we will always assume that the manifold $N$ is a closed
real hyperbolic manifold which is a $\pi$-manifold.
To prove Theorem 3 we will first reduce the problem to another problem.
\vspace{.3in}

\noindent {\bf I. First reduction.}\\

Recall that for any manifold $L$, $DIFF(L)$ denotes the space of all self-diffeomorphisms of $L$,
with the smooth topology,  and if
$\p L\neq\emptyset$,
$DIFF(L,\p )$ denotes the space of all self-diffeomorphisms of $L$
which are the identity on the boundary. \\

Let $\Omega\, DIFF(\D^{n-1},\p )$ be the space of all (continuous) loops in $DIFF(\D^{n-1}, \p )$
based at the identity $1_{\D^{n-1}}$. 
A loop $t\mapsto f_t$, $t\in [0,1]$ is smooth if the map
$(x,t)\mapsto f_t(x)$ is smooth. Classical approximation methods (e.g. using convolution) show that
the inclusion of the space of all smooth loops into
$\Omega\,  DIFF(\D^{n-1},\p )$ is a homotopy equivalence. Hence we will assume, when necessary, that
loops are smooth. We will also assume, if necessary that the loops are constant near 0,1. This does not cause
any problems either.\\

We define a map $$\alpha : \Omega\,  DIFF(\D^{n-1},\p )\ra DIFF(\D^{n},\p )$$

\noindent by the formula $\alpha (f_t)\,\, (x,t)\, =\, (f_t(x),t)$, for $(x,t)\in \D^{n-1}\times [0,1]=\D^n$.
Here $t\mapsto f_t$ denotes a loop in  $DIFF(\D^{n-1})$ and we are identifying $\D^{n-1}\times [0,1]$ with $\D^n$.
(Certainly here we have to assume that the loops are smooth. We also have to smooth corners.)\\

\noindent {\bf Remark.} This map $\alpha$ and the standard constructions of it used here have appeared
(much earlier) in Gromoll's fundamental work \cite{Gm} on positive curvature questions.\\

Identify $\D^{n-1}$ with, say, the northern hemisphere of the sphere $\bS^{n-1}$. Then we have inclusions
$\D^n=\D^{n-1}\times [1,2] \hookrightarrow \bS^{n-1}\times [1,2]\hookrightarrow N$. The composition induces
a map $DIFF(\D^n, \p )\hookrightarrow DIFF_0(N)$, and this map factors through $DIFF_0(\bS^{n-1}\times [1,2],\p )$:
$$DIFF(\D^n, \p )\hookrightarrow DIFF_0(\bS^{n-1}\times [1,2],\p )\stackrel{\Lambda}{\hookrightarrow} DIFF_0(N)$$

\noindent and we denote this composition also by $\Lambda$.\\

\noindent {\bf Remark.} 
As in the remark before the statement of Theorem 3, we will assume that the elements in 
$DIFF(\D^n, \p )=DIFF(\D^{n-1}\times [0,1],\p )$ are constant near $\p ( \D^{n-1}\times [0,1])=
\D^{n-1}\times \{ 0,1\}\cup\D^{n-1}\times [0,1]$. We make this assumption so that the
map $DIFF(\D^n, \p )\hookrightarrow DIFF_0(\bS^{n-1}\times [1,2],\p )$ is 
well defined. Again, as before, this does not cause problems since standard extension methods (along collars) show that the 
inclusion of the space of all elements in $DIFF(\D^n, \p )$ with the properties above into $DIFF(\D^n, \p )$
is a homotopy equivalence.\\

Now, note the following simple but important fact: an element in $DIFF(\D^n,\p )$ is mapped to $\cG$ by the map
$DIFF(\D^n,\p )\ra DIFF_0(\bS^{n-1}\times [1,2],\p )$ if and only if it is in the image of $\alpha$. Therefore we have
reduced the proof of Theorem 3 to the following statement:\\

\noindent {\bf (3.1)} {\it Fix $k\geq 0$. Then, for sufficiently large $n$, $n+k\equiv 2$ $mod \, 4$, the composition map
$$\pi_k(\, \Omega\,  DIFF(\D^{n-1},\p \, ))\stackrel{\alpha_{\#}}{\ra}\pi_k(DIFF(\D^n,\p ))\stackrel{\Lambda_{\#}}{\ra}
\pi_k(DIFF_0(N^n))$$
is non-zero. Also, if $k=0$, $n\geq 10$ and $\Theta_{n+1}\neq 0$, the composition map is also non-zero}.\\

\noindent Here, as we mentioned at the beginning of this section, $N^n$ is any closed real hyperbolic manifold
of dimension $n$ which is a $\pi$-manifold. Note that the last statement of (3.1) corresponds to the Addendum to 
Theorem 3.
We now further reduce statement (3.1) above to another statement.
\vspace{.3in}

\noindent {\bf II. Second reduction.}\\

First we recall some definitions.
For a compact smooth manifold $L$ let $\td (L)$ denote the semi-simplicial group whose $l$-simplices are
self-diffeomorphisms of $\Delta^l\times L$ that send faces $\sigma\times L$ to themselves. Here
$\Delta^l$ is the $l$-simplex and $\sigma$ is any sub-simplex of $\Delta^l$ (see \cite{W}, Sec. 17A). We can consider
$DIFF(L)$ as contained in $\td (L)$ in two (homotopy equivalent) ways: as the set of vertices
of $\td (L)$ or as the semi-simplicial subgroup whose $i$-simplices are
self-diffeomorphisms of $\Delta^l\times L$ that commute with the projection to $\Delta^l$.
Also, define $\td (L,\p )$ as before, but with the extra requirement that
the self-diffeomorphisms of $\Delta^l\times L$ be the identity on $\Delta^l\times \p L$.\\

If we replace ``diffeomorphism'' above by ``homeomorphism'' or ``simple homotopy equivalence''
we obtain spaces ${\widetilde{TOP}}(L)$ and ${\widetilde{G}}(L)$ (and also  ${\widetilde{TOP}}(L,\p )$). 
Here $G(L)$ is the groupoid
of all simple homotopy equivalences of $L$. Since a self-homotopy equivalence does not
have to be one-to-one we have that $G(L)$ and ${\widetilde{G}}(L)$ are homotopy equivalent.
We have fibrations (see \cite{W}, Sec.17A)
$$\td (L)\ra G(L)\ra G/\td (L)$$

$${\widetilde{TOP}}/\td (L)\ra G/\td (L)\ra G/{\widetilde{TOP}}(L)$$
\vspace{.1in}

It is known that $\pi_i({\widetilde{TOP}}/\td (L))\cong [L\times\D^i,\p ;TOP/O]$, where
$[\, ,\, ]$ denotes ``homotopy classes of maps''. Since $TOP/O$ is an infinite loop space
it defines a (non-reduced) generalized cohomology theory such that $h^{-i}(L)=[L\times\D^i,\p ;TOP/O]$.\\

We now come back to the proof. The map $\Lambda : DIFF(\D^n,\p )\ra DIFF_0(N)$ clearly induces a
semi-simplicial map $ \td (\D^n,\p )\ra \td (N)$.\\

\noindent {\bf Lemma 3.2.} {\it $\pi_i(\td (\D^n,\p ))\ra \pi_i(\td(N))$ is a monomorphism for every $i$, provided
$n\geq 5$.}\\

(Recall that we are assuming that $N$ is any closed real hyperbolic manifold
of dimension $n$ which is also a $\pi$-manifold.)\\

\noindent {\bf Proof of the Lemma.}
We claim that
\begin{enumerate}
\item[{1.}] $\pi_i(\td (\D^n,\p ))\cong\pi_{i+1} (G/\td (\D^n,\p ))$.

\item[{2.}] $\pi_i(\td (N))\cong\pi_{i+1} (G/\td (N))$, when $i\geq 1$, and (for $i=0$)
$\pi_{1} (G/\td (N))$ naturally injects into $\pi_0(\td (N))$.

\item[{3.}] $\pi_i({\widetilde{TOP}}/\td (\D^n,\p ))\cong \pi_i(G/\td (\D^n,\p ))$.

\item[{4.}] $\pi_i({\widetilde{TOP}}/\td (N))\cong \pi_i(G/\td (N))$, provided $n\geq 5$.
\end{enumerate}

To prove 1 and 2 use the first fibration above and
just note that  (a) every homotopy equivalence on the disc $\D^n$, mod boundary, can be canonically deformed
(by Alexander's Trick)
to the identity, hence $G(\D^n, \p )$ is contractible; (b) since
$N$ is negatively curved, we have that $G(N)\cong Out(\, \pi_1(N)\,)$, which is a discrete set.
To prove 3 and 4 use the second fibration above and
just note that  (c) $G/{\widetilde{TOP}}(\D^n,\p )$ is  contractible, because the canonical deformation mentioned
above preserves homeomorphisms.
(d) $G/{\widetilde{TOP}}(N)$ is  contractible by Farrell-Jones Rigidity results \cite{FJ2}.\\

It follows from 1-4 above that: $$\pi_i(\td (\D^n,\p ))\cong \pi_{i+1}({\widetilde{TOP}}/\td (\D^n,\p ))\cong [\D^{n+i+1},\p ;TOP/O]
={\tilde{h}}^{-(i+1)}(\bS^n)$$
and $$\pi_i(\td (N ))\cong \pi_{i+1}({\widetilde{TOP}}/\td (N ))\cong [N\times\D^{i+1},\p ;TOP/O]=h^{-(i+1)}(N)$$

\noindent provided $n\geq 5$ and $i\geq 1$.\\

We have then the following commutative diagram:

$$\begin{array}{ccccccc} {\tilde{h}}^{-(i+1)}(\bS^n)&=&\Bigl[\D^{n+i+1},\p ;TOP/O \Bigr]   & \ra & \Bigl[N\times\D^{i+1},\p ;TOP/O\Bigr]&
=&h^{-(i+1)}(N)
\\ &&\downarrow&&\downarrow&&\\  && \pi_i(\td (\D^n,\p ))&\ra& \pi_i(\td(N))&&
\end{array}$$\\

Here the vertical arrows are the canonical isomorphisms mentioned above, and, when $i=0$, the second vertical
is still injective (by 2 above).
Now, since every diffeomorphism in the
image of the map $DIFF (\D^n,\p )\ra DIFF(N)$ is the identity outside $\D^n\sbs N$ we have that the image of an element
in $[\D^{n+i+1},\p ;TOP/O]$ by the map  $[\D^{n+i+1},\p ;TOP/O]\ra [N\times\D^{i+1},\p ;TOP/O]$  in the diagram above
has the property that it is constant outside $\D^n\sbs N$. Hence the map $[\D^{n+i+1},\p ;TOP/O]\ra [N\times\D^{i+1},\p ;TOP/O]$
is induced by the map $c\times 1_{\D^{i+1}}:N\times\D^{i+1}\ra \bS^n\times\D^{i+1}$, where $c:N\ra N/(N\sm\D^n)=\D^n/\p=\bS^n$
is the collapsing map. Therefore, at the cohomology level, the map $h^{-(i+1)}(\bS^n)\ra h^{-(i+1)}(N)$ induced by the
degree-one collapsing map $c$ composed with the canonical monomorphism ${\tilde{h}}^{-(i+1)}(\bS^{n})\ra h^{-(i+1)}(\bS^{n})$
is the homomorphism of Lemma 3.2. The Lemma now follows from the following result: \\

{\it If $N^n$ is a closed stably parallelizable manifold and $f:N\ra\bS^n$ is a degree-one map, then $f^*:h^*(\bS^n)\ra
h^*(N)$ is a monomorphism, for any representable generalized cohomology theory $h$.}\\

For the proof of this result see \cite{FJ1}, claim 2.4, p.902. (Replace $TOP/O$ in the proof of claim
2.4 of \cite{FJ1} by the infinite loop space corresponding to $h$.)
Also, the referee points out that this type of result has long been known to be
an easy corollary to results of G. W. Whitehead \cite{Wh} and W. Browder \cite{Br}. The whole point being
that $\pi$-manifolds are orientable for such theories.
This completes the proof of the Lemma.\\

Note that $\pi_k(\td (\D^n,\p ))\cong \Bigl[\D^{n+k+1},\p ;TOP/O \Bigr]\cong \Theta_{n+k+1}$, the group of homotopy
spheres of dimension $n+k+1$.\\

Consider now the following commutative diagram:
$$\begin{array}{ccccc}\pi_k(\Omega\,  DIFF(\D^{n-1},\p))&\stackrel{\alpha_{\#}}{\ra}&\pi_k(DIFF(\D^n,\p ))&
\stackrel{\Lambda_{\#}}{\ra}&\pi_k(DIFF(N))  \\
& & \downarrow & & \downarrow \\
& &  \pi_k(\td (\D^n,\p ))&\ra& \pi_k(\td(N))
\end{array}$$

\noindent and consider now just the left part of the diagram above:

\begin{equation}\begin{array}{ccc}\pi_k(\Omega\,  DIFF(\D^{n-1},\p))&\stackrel{\alpha_{\#}}{\ra}&\pi_k(DIFF(\D^n,\p ))  \\
& & \downarrow \\
& &  \pi_k(\td (\D^n,\p ))
\end{array}\label{1}\end{equation}\\

From Lemma 3.2 above we see that (3.1) is implied by the following statement:\\

\noindent {\bf (3.3)} {\it Fix $k\geq 0$. Then, for sufficiently large $n$,  $n+k\equiv 2$ $mod \, 4$, the composition map 
$$\pi_k(\Omega\,  DIFF(\D^{n-1},\p))\stackrel{\alpha_{\#}}{\ra}\pi_k(DIFF(\D^n,\p ))\ra \pi_k(\td (\D^n,\p ))$$

\noindent in diagram (1) above is non-zero. Also, if $k=0$, $n\geq 6$ and $\Theta_{n+1}\neq 0$, the composition map is also non-zero.}\\

Note that we have succeeded in eliminating the manifold $N$ from the problem.
\vspace{.3in}

\noindent {\bf III. Proof of statement 3.3.}\\

First we recall some definitions and introduce some notation. For a manifold $L$, the space of smooth
pseudoisotopies of $L$ is denoted by $P(L)$, that is, $P(L)$ consists of all self-diffeomorphisms
of $L\times I$ that are the identity on $L\times\{ 0 \}\cup\p L\times I$. Here $I=[0,1]$.
(The condition of the pseudoisotopies being the identity on $\p L\times I$ is useful but superfluous:
the space of pseudoisotopies of $L$  that are the identity just on $L\times\{ 0 \}$ can be identified
with $P(L)$ by ``bending around corners'', see \cite{Hat}.) Note that $P(L)$ is a group with the composition.
We have stabilization maps $\Sigma :P(L)\ra P(L\times I)$. The direct limit of the sequence 
$P(L)\ra P(L\times I)\ra P(L\times I^2)\ra\dots\,\, $ is called the space of stable pseudoisotopies of $L$, 
and it is denoted by $\cP (L)$. We mention two important facts: (1) $\cP (-)$ is a homotopy functor. Therefore
we get, for example, that $\cP (\D^n)=\cP (*)$, where $*$ denotes a point. (2) the map $\pi_k(P(L))\ra\pi_k(\cP (L))$
is an isomorphism for $k<< dim\, L$, see \cite{I}.\\

Let $\cA_{n-1}$ denote the subgroup of $DIFF(\D^{n-1},\p )$ consisting of all self-diffeomorphisms 
of $(\D^{n-1},\p )$ which are pseudoisotopic
to the identity. Hence we have a surjective homomorphism $\tau : P(\D^{n-1})\ra\cA_{n-1}$, where we just take the ``top'' of a
pseudoisotopy $f$, that is, $\tau (f)=f|_{\D^{n-1}\times\{ 1\}}$. The kernel of this map is the space of all
pseudoisotopies that are the identity on $\p (\D^{n-1}\times I)$, and we can identify this space with
$DIFF(\D^{n},\p )$. Consequently, we get the following sequence:
\vspace{.1in}

\begin{equation} DIFF(\D^{n},\p )\stackrel{\eta}{\ra} P(D^{n-1})\stackrel{\tau}{\ra}\cA_{n-1}
\label{2}\end{equation}
\vspace{.1in}

\noindent {\bf Claim 3.4} {\it  Sequence (2) is a Hurewicz fibration.}\\

Let $t\mapsto f_t$ be a (smooth) path in $\cA_{n-1}$ beginning at $f_0$ and ending at $f_1$. 
(Assume also that the path is constant near $0,1$.) Let
$F_0$ be a lifting of $f_0$, that is, $\tau (F_0)=f_0$. To define a lifting $t\mapsto F_t$ of the whole path
$t\mapsto f_t$, beginning at $F_0$, just take the concatenation of $F_0$ with  the map $(x,t')\mapsto (f_{t'}(x), t')$,
for $0\leq t'\leq t$. (This defines $F_t$ on $\D^{n-1}\times [0,1+t]$ so we need to rescale back to the
interval [0,1].) This proves
the claim.\\

\noindent {\bf Claim 3.5} {\it  The connecting map $\beta : \Omega\, (\cA_{n-1})\ra DIFF(\D^n,\p )$ of 
fibration (2) is homotopic to $\alpha$.}\\

Let $t\mapsto f_t$ be a (smooth) path in $\cA_{n-1}$ beginning at $f_0=1_{\D^{n-1}}$ and ending at $f_1=1_{\D^{n-1}}$. Take
$F_0=1_{\D^n}$ as the lifting of $f_0$. Using the lifting given in the proof of claim 3.4, we see that
the path $t\mapsto f_t$ maps, by the connecting map,  to $F_1$, where $F_1$ is such that $F_1 (x,t)=(x,t)$ for
$0\leq t\leq 1/2$ and $F_1 (x,t)=(f_{2t-1}(x),t)$ for $1/2\leq t\leq 1$. By squeezing the interval $[0,1/2]$ to 0, 
we see that the connecting map $\beta$ is homotopic to $\alpha$. This proves the claim.\\

Now, to prove (3.3) it is enough to prove the following:\\

\noindent {\bf (3.6)} {\it Fix $k\geq 0$. Then, for sufficiently large $n$, $n+k\equiv 2$ $mod \,4$, the composition map 
$$\pi_k(\Omega\, \cA_{n-1})\stackrel{\beta_{\#}}{\ra}\pi_k(DIFF(\D^n,\p ))\ra \pi_k(\td (\D^n,\p ))$$

\noindent  is non-zero. Also, if $k=0$, $n\geq 10$ and $\Theta_{n+1}\neq 0$, the composition map is also non-zero.}\\

Using the fibration   $DIFF(\D^n,\p )\ra\td (D^n, \p )\ra \frac{\td}{DIFF}(\D^n,\p )$ (see \cite{W}, Sec. 17A) and fibration (2),
we can embed the sequence in (3.6) into the following larger diagram:

\begin{equation}\begin{array}{ccccc}\pi_k(\Omega\,  \cA_{n-1})&\stackrel{\beta_{\#}}{\ra}&\pi_k(DIFF(\D^n,\p ))&
\stackrel{}{\ra}&\pi_k(P(\D^{n-1}))  \\
& & \downarrow & &  \\
& &  \pi_k(\td (\D^n,\p ))&=& \Theta_{n+k+1}\\
& & \downarrow & &  \\
& &  \pi_k(\frac{\td}{DIFF} (\D^n,\p ))&& 
\end{array}\label{3}\end{equation}
\vspace{.1in}

\noindent where the upper row and the central column are exact. 
Note that for $n>>k$, $\pi_k(P(\D^{n-1}))\cong  \pi_k(\cP (\D^{n-1}))\cong \pi_k(\cP (*))$ which does not
depend on $n$ (see \cite{I}).
To prove (3.6)
we have three cases, with increasing degrees of difficulty: the case $k=0$ (which corresponds to
the Addendum to Theorem 3), the case $k\not\equiv 3$ $mod\, 4$, and the case $k\equiv 3$ $mod\, 4$.
\vspace{.4in}

\noindent {\bf First case: $k=0$, $n\geq 6$ and $\Theta_{n+1}\neq 0$.}\\

In this case, since $n-1\geq 5$ we have that $\pi_0(P(\D^{n-1}))\cong 0$ by Cerf's foundational work \cite{Cf}. 
Also, it follows immediately from the definitions that $\pi_0(\frac{\td}{DIFF} (\D^n,\p ))\cong 0$.
Since we are assuming that $\Theta_{n+1}\neq 0$, we can (using diagram (3)) pull-back a non-zero element $a\neq 0\in \Theta_{n+1}$
all the way back to $\pi_0(\Omega\,  \cA_{n-1})$. This proves (3.6) for this first case.
\vspace{.4in}

\noindent {\bf Second case: $k\not\equiv 3$ $mod\, 4$, $n>>k$ and $n+k\equiv 2$ $mod\, 4$.}\\

We will use the following facts:
\begin{enumerate}
\item[{1.}] For $k\not\equiv 3$ $mod\, 4$, $\pi_k(\cP(*))$ is a finite group. For 
$k\equiv 3$ $mod\, 4$, $\pi_k(\cP(*))\cong\Z\oplus$(finite group$)$ (see \cite{Dw}, \cite{FH}, \cite{Wa}). 
Denote by $a_k$ the order of the torsion part of $\pi_k(\cP(*))$.

\item[{2.}] Using Hatcher's spectral sequence (see \cite{Hat} prop. 2.1, 2.2) we have that, for 
$n>>k$  the group $\pi_k(\frac{\td}{DIFF} (\D^n,\p ))$ has a filtration $0=G_0<...<G_k=\pi_k(\frac{\td}{DIFF} (\D^n,\p )$
such that $G_{i}/G_{i-1}$ is a subquotient of $H_{k-i}(\Z_2,\pi_{i-1}\cP (*)))$. Since all groups $H_{k-i}(\Z_2,\pi_{i-1}\cP (*)))$
are 2-torsion groups, for $i<k$ and $H_{0}(\Z_2,\pi_{k-1}\cP (*))\cong \pi_{k-1}\cP (*)/\Z_2$
(the quotient by the action of $\Z_2$), we see that 
the torsion part of $\pi_k(\frac{\td}{DIFF} (\D^n,\p ))$ has  $2^{k}a_{k-1}$ for an exponent.

\item[{3.}] The group $\Theta_{4m-1}$ has a cyclic subgroup of order $2^{2m-2}(2^{2m-1}-1)numerator\{ 4B_{m}/m\}$
(see \cite{M}, p.285). Note that this order increases exponentially with $m$.
\end{enumerate}

The important observation here is that, for $n>>k$, the group $\pi_k(P(\D^{n-1}))$ and an exponent of the torsion part of 
$\pi_k(\frac{\td}{DIFF} (\D^n,\p ))$ do not depend on $n$.\\

\noindent {\bf Remarks.}

\noindent {\bf 1.} Here and in the next case (that is, in the proofs of the second and third cases) item 2 above can be
replaced by the following fact that can be deduced using \cite{BL1} or Lemma 2.2 of \cite{HJ}:
{\it the exponent of the odd order torsion part of $\pi_k(\frac{\td}{DIFF} (\D^n,\p ))$, $n>>k$, does not depend
on $n$}. (Indeed this exponent is $a_{k-1}$, where $a_{k-1}$ is as in item 1.) To use this fact instead of
item 2 just note that by item 3, $\Theta_{4m-1}$ has elements of order $2^{2m}-1$,
which is large and odd.\\

\noindent {\bf 2.} Antonelli, Burghelea and Kahn \cite{ABK} showed that the image of
$$\pi_k(\, DIFF (\D^n, \p ) \, )\ra  \pi_k(\frac{\td}{DIFF} (\D^n,\p ))$$
is isomorphic to Gromoll's group $\Gamma_{k+1}^{n+k+1}$ (see \cite{Gm}).
Furthermore they obtained many strong non-vanishing results in \cite{ABK} about these
groups of Gromoll. Their results perhaps combined with recent knowledge about $a_k$
(see \cite{Rg1}, \cite{Rg2}) should yield extremely substantial quantitative improvements to
Theorem 3 and hence also to Theorem 1. We are very grateful to the referee for pointing out
to us this direction for future investigation.\\

We continue with the proof of the second case.
Fix $k$ with $k\not\equiv 3$ $mod\, 4$.
Now, since we are assuming $n+k\equiv 2$ $mod\, 4$, we have that $n+k+1=4m-1$, for some $m$.
Hence, using diagram (3) and the three facts above, we get that, by choosing $n$ large,  we can find a non-zero element $x\in
\pi_k(\td (\D^n,\p ))= \Theta_{n+k+1}=\Theta_{4m-1}$ with large order that maps to $0\in \pi_k(\frac{\td}{DIFF} (\D^n,\p ))$.
Hence $x$ is the image of an element $y\in\pi_k(DIFF(\D^n,\p ))$ with large order. Now, since we are assuming
$k\not\equiv 3$ $mod\, 4$, by fact 1 above we have that $\pi_k(P(\D^{n-1}))\cong \pi_k(\cP (*))$ is a finite group
of order $a_k$. Then $a_ky$ maps to $0\in\pi_k(P(\D^{n-1}))$, hence $a_ky$ pulls-back to an element $z$
in $\pi_k(\Omega\,  \cA_{n-1})$.
Since we can take the order of $x$ and $y$ as large as we
want, we can choose $x$ and $y$ such that $a_k x\neq 0$. Hence $a_ky\neq 0$ and it follows that $z\neq 0$. This concludes
the proof of (3.6) in the second case.\\

\noindent {\bf Third case: $k\equiv 3$ $mod\, 4$, $n>>k$ and $n+k\equiv 2$ $mod\, 4$.}\\

The problem in this case is that now the group $\pi_k(P(\D^{n-1}))\cong\Z\oplus$(finite group$)$, hence it
is not finite, and we cannot use the argument above because the element  $y$ can map to an infinite order element.\\

To begin with, we embed diagram (3) in a larger diagram:

\begin{equation}\begin{array}{ccccc}\pi_k(P(\D^{n-1}))\\ \Sigma_{\#}\downarrow \\
\pi_k(P(\D^{n}))&\stackrel{\phi}{\ra}&\pi_{k+1}(\frac{\td}{DIFF} (\D^n,\p ))\\
 &  \stackrel{\tau_{\#}}{\searrow} & \mu\downarrow\\
\pi_k(\Omega\,  \cA_{n-1})&\stackrel{}{\ra}&\pi_k(DIFF(\D^n,\p ))&
\stackrel{\eta_{\#}}{\ra}&\pi_k(P(\D^{n-1}))  \\
& & \nu\downarrow & &  \\
& &  \pi_k(\td (\D^n,\p ))&=& \Theta_{n+k+1}\\
& & \downarrow & &  \\
& &  \pi_k(\frac{\td}{DIFF} (\D^n,\p ))&& 
\end{array}\label{4}\end{equation}
\vspace{.1in}

We explain the new terms. The central column contains one more term than the central
column of (3) and, as before, it is a piece of the exact sequence of the homotopy groups of the fibration
$DIFF(\D^n,\p )\ra\td (D^n, \p )\ra \frac{\td}{DIFF}(\D^n,\p )$.
The map $\Sigma_{\#}$ is induced by the suspension map $\Sigma :P(\D^{n-1})\ra P(\D^{n-1}\times I)=P(\D^n)$.
The diagonal arrow is induced the map $\tau: P(D^n)\ra DIFF(\D^n)$, which, as defined before, consists
of taking the ``top'' of the pseudoisotopy. An element $f$ in $\pi_{k+1}(\frac{\td}{DIFF} (\D^n,\p ))
=\pi_{k+1}(\td (\D^n,\p ),DIFF(\D^n,\p ))$ is represented by a self-diffeomorphism of $\D^n\times\D^{k+1}$ which
preserves the projection to $\D^{k+1}$ over $\p\, \D^{k+1}$ and is the identity when restricted to
$(\p\D^n\times \D^{k+1})\cup (\D^n\times\{ 0\}\times\D^k)$, where $\D^{k+1}=I\times\D^k$. 
Hence $\mu(f)$ is represented by the restriction
of the above self-diffeomorphism to $\D^n\times\p\, \D^{k+1}$. An element in $\pi_k(P(\D^{n}))$
is represented by a self-diffeomorphism of $(\D^n\times I)\times\D^k$ which is the identity
over $\p\Bigl( (\D^n\times I)\times\D^k\Bigr)\sm (\D^n\times \{ 1\}\times\D^k)$, and preserves the
projection to $\D^k$. Identifying $I\times \D^k$ with $\D^{k+1}$ we obtain a map 
$\pi_k(P(\D^{n}))\ra\pi_{k+1}(\frac{\td}{DIFF} (\D^n,\p ))$. This map is $\phi$. It is easy to verify that
$\mu\phi=\tau_{\#}$. Hence diagram (4) is commutative with the central row and central column exact.\\

Recall that we are assuming that $k\equiv 3$ $mod\, 4$, and that $n+k\equiv 2$ $mod\, 4$.
It follows that $n\equiv 3$ $mod \, 4$. In particular $n$ is odd. Then,
assuming also that $n>>k$, we have the following facts:
\begin{enumerate}
\item[{i.}] $\Sigma_{\#}$ is an isomorphism (see\cite{I})
and $\pi_k(P(\D^n))\cong \pi_k(P(\D^{n-1}))\cong\pi_k(\cP (*))\cong\Z\oplus($finite group$)$
(see \cite{Dw}, \cite{FH}, \cite{Wa}). An element in this group can then be written in the form $j+t$, where $j\in \Z$ and
$t$ is in the torsion part of the group.

\item[{ii.}] $\eta_{\#}\, \tau_{\#}\Sigma_{\#}(x)=x\pm\bar x$, for $x\in \pi_k(P(\D^{n-1}))$. Here the map
$x\mapsto\bar x$ is an involution on $\pi_k(P(\D^{n-1}))$ (see \cite{Hat}).

\item[{iii.}] $\pi_k(DIFF(\D^n,\p ))\otimes\Q\cong\Q$ and $\tau_{\#}\otimes \Q$ 
is an isomorphism (see \cite{FH}, Th. 2.1). Also $\eta_{\#}\otimes \Q$ is an isomorphism.
(For this use also Th. 2.1 of \cite{FH} and combine it with the homotopy exact sequence induced
by the fibration (2). One detail: note that $\pi_i(\cA_{n-1})\cong\pi_i(DIFF(\D^n,\p ))$, for $i>0$,
because $DIFF^0(\D^n,\p )\sbs \cA_{n-1}\sbs DIFF(\D^n,\p )$, where $DIFF^0(\D^n,\p )$ is the
connected component of the identity in $DIFF(\D^n,\p )$.)
\end{enumerate}

Define $y_0=\tau_{\#}\Sigma_{\#}(1)\in\pi_k(DIFF(\D^n,\p ))$. Then, by (i) and (iii) above,
$y_0$ has infinite order, and, since diagram (4) is commutative, $\nu(y_0)=0$.
We have $\eta_{\#}(y_0)=\eta_{\#}\, \tau_{\#}\Sigma_{\#}(1)$, and by (ii), 
$\eta_{\#}(y_0)=1\pm\bar 1$. But an involution sends 1 to $\pm 1+t$, where
$t$ is a torsion element. Hence $\eta_{\#} (y_0)$ is either of the form $2+t$ or $t$, and, since $y_0$ has infinite
order, (iii) shows that $\eta_{\#} (y_0)=2+t$.
Write $y_1=a_ky_0$. Then $\eta_{\#} (y_1)=2a_k$. Note that
$y_1$ also has infinite order and $\nu(y_1)=0$.
As in the second case (i.e for $k\not\equiv 3$ $mod\, 4$), we can find an element $y'\in \pi_k(\, DIFF(\D^n,\p )\, )$
such that $\nu(y')\neq 0$ and  $\eta_{\#}(y')=2a_kj$ ($y'=2a_ky$, where $y$ is as in the proof of the second case).
Now take $y''=y'-jy_1$ and we see that $\eta_{\#}(y'')=0$ and $\nu(y'')=\nu(y')\neq 0$, and we are done.
This completes the proof of Theorem 3.

F.T. Farrell

SUNY, Binghamton, N.Y., 13902, U.S.A.\\

P. Ontaneda

SUNY, Binghamton, N.Y., 13902, U.S.A.

\end{document}